\documentclass[11pt,reqno,tbtags]{article}
\usepackage{amsmath}
\usepackage{amsthm}

\numberwithin{equation}{section}

\newtheorem{theorem}{Theorem}[section]
\newtheorem{corollary}[theorem]{Corollary}
\newtheorem{lemma}[theorem]{Lemma}
\newtheorem{proposition}[theorem]{Proposition}
\newtheorem{claim}[theorem]{Claim}
\newtheorem{example}[theorem]{\sl Example}

\theoremstyle{definition}
\newtheorem{remark}[theorem]{Remark}
\newtheorem{problem}[theorem]{Open Problem}

\newcommand{\lc}{\left\lceil}
\newcommand{\lf}{\left\lfloor}
\newcommand{\rc}{\right\rceil}
\newcommand{\rf}{\right\rfloor}

\newcommand{\EE}{{\bf  E}}

\newcommand{\RR}{{\bf  R}}

\newcommand{\TT}{{\bf  T}}
\newcommand{\ZZ}{{\bf  Z}}

\newcommand{\Var}{{\bf Var}}

\newcommand{\Ct}{{\tilde{C}}}

\newcommand{\Yt}{{\tilde{Y}}}

\newcommand{\Fc}{{\cal F}}

\newcommand{\Lc}{{\cal L}}
\newcommand{\Lto}{{\stackrel{\Lc}{\to}}}
\newcommand{\Leq}{{\,\stackrel{\Lc}{=}\,}}
\newcommand{\Lnoteq}{{\,\stackrel{\Lc}{\ne}\,}}

\newcommand{\Yh}{\hat{Y}}

\newcommand{\nuh}{\hat{\nu}}

\newcommand{\begp}{\begin{proposition}}
\newcommand{\enp}{\end{proposition}}
\newcommand{\begt}{\begin{theorem}}
\newcommand{\ent}{\end{theorem}}
\newcommand{\begl}{\begin{lemma}}
\newcommand{\enl}{\end{lemma}}
\newcommand{\begc}{\begin{corollary}}
\newcommand{\enc}{\end{corollary}}
\newcommand{\begcl}{\begin{claim}}
\newcommand{\encl}{\end{claim}}
\newcommand{\begr}{\begin{remark}}
\newcommand{\enr}{\end{remark}}
\newcommand{\begal}{\begin{algorithm}}
\newcommand{\enal}{\end{algorithm}}
\newcommand{\begd}{\begin{definition}}
\newcommand{\enf}{\end{definition}}
\newcommand{\begx}{\begin{example}}
\newcommand{\enx}{\end{example}}
\newcommand{\bega}{\begin{array}}
\newcommand{\ena}{\end{array}}

\newcommand{\sfrac}[2]{{\textstyle\frac{#1}{#2}}}

\def\rompar(#1){\textup(#1\textup)}    

\newcommand\eps{\varepsilon}

\newcommand\gl{\lambda}

\newcommand\abar{\bar{a}}
\newcommand\Mbar{\bar{M}}

\newcommand\Vbar{\bar{V}}

\newcommand\fudge{[1 + (1 / n)]}

\newcommand{\refS}[1]{Section~\ref{#1}}
\newcommand{\refT}[1]{Theorem~\ref{#1}}
\newcommand{\refC}[1]{Corollary~\ref{#1}}
\newcommand{\refL}[1]{Lemma~\ref{#1}}
\newcommand{\refR}[1]{Remark~\ref{#1}}

\renewcommand\Re{\operatorname{Re}}

\newcommand\nopf{\qed}   
\newcommand\noqed{\renewcommand{\qed}{}} 
\newcommand\qedtag{\tag*{\qedsymbol}}

\newcommand\CS{Cauchy--Schwarz}
\newcommand\KSm{Kolmogorov--Smirnov}

\newcommand\Roesler{R\"{o}sler}
\newcommand\Quicksort{\texttt{Quicksort}}

\newcommand\dks{d_{\mbox{\rm \scriptsize KS}}}
\begin{document}

\setcounter{page}{0}
\thispagestyle{empty}

\begin{center}
{\Large \bf Quicksort Asymptotics \\}
\normalsize

\vspace{4ex}
{\sc James Allen Fill\footnotemark} \\
\vspace{.1in}
Department of Mathematical Sciences \\
\vspace{.1in}
The Johns Hopkins University \\
\vspace{.1in}
{\tt jimfill@jhu.edu} and {\tt http://www.mts.jhu.edu/\~{}fill/} \\
\vspace{.2in}
{\sc and} \\
{\sc Svante Janson}\\ 
\vspace{.1in}
Department of Mathematics \\
\vspace{.1in}
Uppsala University \\
\vspace{.1in}
{\tt svante.janson@math.uu.se} and {\tt http://www.math.uu.se/\~{}svante/} \\
\end{center}
\vspace{3ex}

\begin{center}
{\sl ABSTRACT} \\
\end{center}

The number of comparisons~$X_n$ used by \Quicksort{} to sort an array
of~$n$ distinct numbers has mean $\mu_n$ of order $n \log n$ and standard
deviation of order~$n$.  Using different methods, R\'{e}gnier and \Roesler{}
each showed that the normalized variate $Y_n := (X_n - \mu_n) / n$ converges in
distribution, say to~$Y$; the distribution of~$Y$ can be characterized as the
unique fixed point with zero mean of a certain distributional transformation.

We provide the first rates of convergence for the distribution of~$Y_n$ to
that of~$Y$, using various metrics.  In particular, we establish the bound
$2 n^{- 1 / 2}$ in the $d_2$-metric, and the rate $O(n^{\eps - (1 / 2)})$ for
\KSm{} distance, for any positive~$\eps$.
\bigskip
\bigskip

\begin{small}

\par\noindent
{\em AMS\/} 2000 {\em subject classifications.\/}  Primary 68W40;
secondary 68P10, 60F05, 60E05.
%
\medskip
\par\noindent
\emph{Date.} May~22, 2001. 
\end{small}

\footnotetext[1]{Research supported by NSF grant DMS--9803780,
and by The Johns Hopkins University's Acheson J.~Duncan Fund for the
Advancement of Research in Statistics.}

\newpage
\addtolength{\topmargin}{+0.5in}

\section{Introduction and summary}
\label{S:intro}

This  
paper provides the first rates of convergence (as $n \to \infty$) for the
distribution of the number of comparisons used by the sorting algorithm
\Quicksort{} to sort an array of~$n$ distinct numbers.  \Quicksort{} is the
standard sorting procedure in {\tt Unix} systems, and has been
cited~\cite{DS} as one of the ten algorithms
``with the greatest influence on the development and practice
of science and engineering in the 20th century."  We begin with a brief review
of what is known about the analysis of {\tt Quicksort} and a summary of our
new results.

The {\tt Quicksort} algorithm for sorting an array of~$n$ distinct
numbers is
extremely simple
to describe.  If $n = 0$ or $n = 1$, there is nothing to do.  If $n
\geq 2$, pick a number
uniformly at random from the given array.  Compare the other numbers
to it to partition the
remaining numbers into two subarrays.  Then recursively invoke {\tt
Quicksort} on each of
the two subarrays.

Let~$X_n$ denote the (random) number of comparisons required (so that
$X_0 = 0$).  Then~$X_n$
satisfies the distributional recurrence relation
\begin{equation}
\label{zeus}
X_n \Leq X_{U_n - 1} + X^*_{n - U_n} + n - 1,\qquad n \geq 1,
\end{equation}
where~$\Leq$ denotes equality in law (i.e.,\ in distribution), and
where, on the right,
$U_n$ is distributed uniformly on the set $\{1, \ldots, n\}$,
$X_j^* \Leq X_j$,
and
$$
U_n;\ X_0, \ldots, X_{n - 1};\ X^*_0, \ldots, X^*_{n - 1}
$$
are all independent.

As is well known and quite easily established, for $n \geq 0$ we have
$$
\mu_n := \EE\,X_n = 2 (n + 1) H_n - 4 n \sim 2 n \ln n,
$$
where $H_n := \sum_{k = 1}^n k^{-1}$ is the $n$th harmonic number
and~$\sim$ denotes
asymptotic equivalence.  It is also routine to compute explicitly the variance
of~$X_n$ (see Exercise 6.2.2-8 in~\cite{Knuth3}):
\begin{equation}
\label{poseidon}
\Var\,X_n = 7 n^2 - 4 (n + 1)^2 H^{(2)}_n - 2 (n + 1) H_n + 13 n = \sigma^2
n^2 - 2 n \ln n + O(n)
\end{equation}
where $H^{(2)}_n := \sum_{k = 1}^n k^{-2}$ are the second-order harmonic
numbers and
\begin{equation}
\label{artemis}
\sigma^2 := 7 - \sfrac{2}{3} {\pi}^2 \doteq 0.42.
\end{equation}

Consider the normalized variate
$$
Y_n := (X_n - \mu_n) / n, \qquad n \geq 1.
$$
Then~\eqref{zeus} implies the recursion
\begin{equation}
\label{athena}
Y_n \Leq \frac{U_n - 1}{n} Y_{U_n - 1} + \frac{n - U_n}{n} Y^*_{n - U_n}
  + C_n(U_n), \qquad n \geq 1,
\end{equation}
with~$Y_0$ arbitrarily defined (since its coefficient is~$0$), where on the
right, as for~$X_n$, we have $U_n \sim \mbox{unif}\{1, \ldots, n\}$ and $Y_j^*
\Leq Y_j$, and $U_n;\ Y_1, \ldots, Y_{n - 1};\ Y^*_1, \ldots, Y^*_{n -
1}$ are all independent; further,
\begin{equation}
\label{hermes}
C_n(i) := \mbox{$\frac{n - 1}{n}$} + \mbox{$\frac{1}{n}$}
(\mu_{i -1} + \mu_{n - i} - \mu_n),\ \ \ 1 \leq i \leq n.
\end{equation}
Note that $\EE\,Y_n = 0 = \EE\,C_n(U_n)$.  We will see below that if $n \to
\infty$ and $i / n \to u \in [0, 1]$, then $C_n(i) \to C(u)$, where 
$$
C(u) := 2 u \ln u + 2 (1 - u) \ln(1 - u) + 1, \qquad u \in [0, 1],
$$
with  the natural (continuous) interpretation $C(u) := 1$ for $u = 0, 1$.

Moreover, R\'{e}gnier~\cite{Reg} and \Roesler~\cite{Roesler} showed, using
different methods, that $Y_n \to Y$ in distribution, with~$Y$ satisfying the
distributional identity
\begin{equation}
\label{fix}
Y\,\Leq\,U Y + (1 - U) Y^* + C(U)
\end{equation}
obtained by formally taking limits in~\eqref{athena}, where,
on the right,
$U$, $Y$, and~$Y^*$ are independent, with $Y^* \Leq Y$  and $U \sim
\mbox{unif}(0, 1)$.  [\Roesler~\cite{Roesler} showed further that~\eqref{fix}
characterizes the limiting law~$\Lc(Y)$, subject to $\EE\,Y = 0$ and $\Var\,Y
< \infty$.  For a complete characterization of the distributions
satisfying~\eqref{fix}, see~\cite{SJFill1A}.]

The purpose of the present paper is to study the rate of convergence of
$\Lc(Y_n)$ to $\Lc(Y)$, using several different measures of the distance
between $\Lc(Y_n)$ and $\Lc(Y)$.

First, for real $1 \leq p < \infty$, let $\| X \|_p := \left( \EE\,X^p
\right)^{1 / p}$ denote the $L^p$-norm, and let $d_p$ denote the metric on the
space of all probability distributions with finite $p$th absolute moment
defined by
$$
d_p(F, G) := \min \| X - Y \|_p,
$$
taking the minimum over all pairs of random variables~$X$ and~$Y$ (defined on
the same probability space) with $\Lc(X) = F$ and $\Lc(Y) = G$.  We will use
the fact~\cite{Cambanis} that the minimum is attained for each $1 \leq p <
\infty$ by the same~$X$ and~$Y$, viz.,\ $X := F^{-1}(u)$ and $Y := G^{-1}(u)$
defined for~$u$ in the probability space $(0, 1)$ (with Lebesgue measure).

We will for simplicity write $d_p(X, Y) := d_p(\Lc(X), \Lc(Y))$ for random
variables~$X$ and~$Y$, but note that this distance depends only on the marginal
distributions of~$X$ and~$Y$.

\Roesler~\cite{Roesler} showed that $d_p(Y_n, Y) \to 0$ as $n \to \infty$
for every $1 \leq p < \infty$.  In Sections~\ref{S:d2} and~\ref{S:dp} we will
quantify this and show that
$$
d_p(Y_n, Y) = O\left(n^{- 1 / 2}\right)
$$
for every fixed~$p$.  In the case $p = 2$ we will further show the explicit
bound
$$
d_2(Y_n, Y) < 2 n^{- 1 / 2}.
$$
We do not know whether the $n^{- 1 / 2 }$ rate is sharp, although that is
widely believed.  The best lower bound we can show
(\refS{S:lower}) is
$$
d_p(Y_n, Y) \geq c \frac{\ln n}{n}, \qquad p \geq 2,
$$
with $c > 0$ independent of~$p$.

In \refS{S:KS} we use these results to bound the \KSm\ distance $\dks(Y_n, Y)$
between $\Lc(Y_n)$ and $\Lc(Y)$.  We show that
$$
\dks(Y_n, Y) = O\left(n^{\eps - (1 / 2)}\right)
$$
for every $\eps > 0$.  We believe that the rate is in fact
$O\left(n^{- 1 / 2}\right)$, but again do not know the exact rate.  The best
lower bound we can prove is $c / n$ with $c > 0$.

In \refS{S:local} we prove a kind of local limit theorem which enables us to
approximate the density function~$f$ of~$Y$.  (It was proved by Tan and
Hadjicostas~\cite{TanH} that~$Y$ has a density function; $f$ is bounded and
infinitely differentiable by~\cite{SJFill1}.)

\Roesler~\cite{Roesler} showed that (for fixed $\gl \in \RR$) the moment
generating function values $\EE\,e^{\gl Y_n}$ are bounded and thus
converge to $\EE\,e^{\gl Y}$.  Again we quantify his bounds and give in
\refS{S:mgf} explicit bounds, based on \Roesler's method.

In several (but not all) bounds we give explicit numerical values to
constants.  These values are hardly the best possible, but we make some effort
to get fairly small values.  This includes sometimes the use of extensive
numerical verifications by computer for small~$n$.  [All  numerical
calculations have been verified independently by the two authors, the
(alphabetically) first using {\tt Mathematica} and the second using
{\tt Maple}.]  Such arguments could be simplified or omitted at the cost of
increasing the constants.

\subsection{Preliminaries} \label{S:prel}

In order to later estimate $C_n(i)$ defined
by~\eqref{hermes} we need some explicit bounds on~$\mu_n$.  First, as
mentioned above,
\begin{equation}
\label{castor}
\mu_n = 2 (n + 1) H_n - 4 n,
\end{equation}
which can be rewritten
\begin{equation}
\label{pollux}
\mu_n = 2 (n + 1) H_{n + 1} - 4n - 2.
\end{equation}
Next we use the bounds on the harmonic numbers (see, e.g.,\ Section~1.2.11.2
in~\cite{Knuth1})
\begin{equation}
\label{jason}
\ln n + \gamma \leq H_n \leq \ln n + \gamma + \sfrac{1}{2 n}, \qquad n \geq 1.
\end{equation}
Hence, for $n \geq 1$, from~\eqref{castor}
\begin{equation}
\label{helen}
2 (n + 1) \ln n + (2 \gamma - 4) n + 2 \gamma \leq \mu_n \leq
2 (n + 1) \ln n + (2 \gamma - 4) n + 2 \gamma + \sfrac{n + 1}{n}
\end{equation}
and from~\eqref{pollux}
\begin{equation}
\label{clymenestra}
2 n \ln n + (2 \gamma - 4) n + 2 \leq \mu_{n - 1} \leq
2 n \ln n + (2 \gamma - 4) n + 3.
\end{equation}

\section{Bounding $d_2(Y_n, Y)$} \label{S:d2}

In this section we prove the following explicit estimate for $d_2(Y_n, Y)$.

\begin{theorem}\label{T:d2}
For all $n \geq 1$, $d_2(Y_n, Y) < 2 / \sqrt{n}$.
\end{theorem}

\begin{proof}
We basically follow the method of \Roesler~\cite{Roesler}, making all
estimates explicit.  We study in this paper only properties of the
univariate distributions of~$Y_n$.  We thus take the liberty of letting~$Y_n$
denote any random variable with the appropriate distribution [$Y_n \Leq (X_n -
\mu_n) / n$].  We then may choose $Y_0, Y_1, \ldots$ defined on the same
probability space as~$Y$ and such that
$$
\| Y_n - Y \|_2 = d_2(Y_n, Y), \qquad n \geq 0.
$$
Further, let $(Y^*, Y_0^*, Y_1^*, \ldots)$ be an independent copy of $(Y, Y_0,
Y_1, \ldots)$ and let $U \sim \mbox{unif}(0, 1)$ be independent of
everything else.  For convenience we write $a_n := d_2(Y_n, Y)$.

Observe, by~\eqref{athena},  that
\begin{equation}
Y_n \Leq \Yt_n := \frac{\lc n U \rc - 1}{n} Y_{\lc n U \rc - 1}
                    + \frac{n - \lc n U \rc}{n} Y^*_{n - \lc n U \rc}
                    + C_n(\lc n U \rc)
\end{equation}
and recall  from~\eqref{fix} that
\begin{equation}
Y \Leq \Yt := U Y + (1 - U) Y^* + C(U).
\end{equation}
Therefore,
\begin{equation}
\label{anbound}
a_n^2 = d^2_2(Y_n, Y) \leq \EE |\Yt_n - \Yt|^2.
\end{equation}
Now
\begin{equation*}
\begin{split}
{\hspace{-2in}}\Yt_n - \Yt
&=
 \left( \frac{\lc n U \rc - 1}{n} Y_{\lc n U \rc - 1} - U Y \right)
   + \left( \frac{n - \lc n U \rc}{n} Y^*_{n - \lc n U \rc}
            - (1 - U) Y^* \right) \\
&\qquad \qquad
   + \left( C_n(\lc n U \rc) - C(U) \right) \\
&  =: W_1 + W_2 + W_3,
\end{split}
\end{equation*}
say.  Given~$U$, the random variables~$W_1$ and~$W_2$ are independent with
zero mean, while~$W_3$ is a constant.  Hence
$$
\EE \left( \left. \left| \Yt_n - \Yt \right|^2\,\right|\,U \right) = \EE
\left( \left. (W_1 + W_2 + W_3)^2 \right|\,U \right) = \EE\left( W^2_1\,|\,
U \right) + \EE\left( W^2_2\,|\,U \right) + W^2_3
$$
and thus, taking expectations,
\begin{equation}
\label{othello}
\EE\left| \Yt_n - \Yt \right|^2 = \EE\,W^2_1 + \EE\,W^2_2 + \EE\,W^2_3.
\end{equation}
By symmetry (replacing~$U$ by $1 - U$), $\EE\,W^2_2 = \EE\,W^2_1$.  We
estimate this term by conditioning on~$U$, using the independence
of~$U$ and $Y, Y_0, \ldots$.  If $U = (k + v) / n$, with $k \in \{0,
1, \ldots, n - 1\}$ and $0 < v \leq 1$, then $\lc n U \rc = k + 1$ and
$W_1 = \sfrac{k}{n} (Y_k - Y) - \sfrac{v}{n} Y$; hence Minkowski's
inequality yields
\begin{eqnarray*}
\EE \left( W^2_1\,|\,U = (k + v) / n \right)^{1 / 2}
 &\leq& \sfrac{k}{n} \| Y_k - Y \|_2 + \sfrac{v}{n} \| Y \|_2 \\
 &  = & \sfrac{k}{n} a_k + \sfrac{v}{n} \sigma.
\end{eqnarray*}
Consequently,
\begin{eqnarray}
\EE\,W^2_1
 &  = & \frac{1}{n} \sum_{k = 0}^{n - 1} \int_0^1\!\EE
          \left(W^2_1\,|\,U = (k + v) / n \right)\,dv
  \leq  \frac{1}{n} \sum_{k = 0}^{n - 1} \int_0^1\!\left( \frac{k}{n}
          a_k + \frac{v}{n} \sigma \right)^2\,dv \nonumber \\
 &  = & \frac{1}{n} \sum_{k = 0}^{n - 1} \int_0^1\!\left(
          \frac{k^2}{n^2} a^2_k + \frac{2 k}{n^2} v a_k \sigma +
          \frac{v^2}{n^2} \sigma^2 \right)\,dv \nonumber \\
\label{desdemona}
 &  = & \frac{1}{n} \sum_{k = 0}^{n - 1} \left( \frac{k^2}{n^2} a^2_k +
          \frac{k}{n^2} a_k \sigma + \frac{\sigma^2}{3 n^2} \right).
\end{eqnarray}
We postpone the estimation of $\EE\,W^2_3$, and introduce the notation
\begin{equation}
\label{iago}
b_n := \| W_3 \|_2 = \| C_n(\lc n U \rc) - C(U) \|_2.
\end{equation}

Combining~\eqref{anbound}--\eqref{iago},
we obtain our fundamental recursive estimate
\begin{eqnarray}
a^2_n &\leq& 2 \EE\,W^2_1 + \EE\,W^2_3 \nonumber \\
\label{shylock}
      &\leq& \frac{2}{n^3} \sum_{k = 1}^{n - 1} k^2 a^2_k + \frac{2
               \sigma}{n^3} \sum_{k = 1}^{n - 1} k a_k + \frac{2 \sigma^2}{3
               n^2} + b^2_n, \qquad n \geq 1.
\end{eqnarray}

We unwrap this recursion partly, by concentrating on the first sum on the
right-hand side and regarding the second as known.  Thus, writing
\begin{equation}
\label{lear}
y_n := \frac{2 \sigma}{n^3} \sum_{k = 1}^{n - 1} k a_k + \frac{2 \sigma^2}{3
n^2} + b^2_n,
\end{equation}
we define recursively
\begin{equation}
\label{portia}
x_n := \frac{2}{n} \sum_{k = 1}^{n - 1} x_k + n^2 y_n, \qquad n \geq 1,
\end{equation}
and find by~\eqref{shylock} and induction
$$
n^2 a^2_n \leq x_n, \qquad n \geq 1.
$$
Now, the recursion~\eqref{portia} is easily
solved (see, e.g.,\ \cite{Fill}),  giving 
\begin{equation}
\label{duncan}
a^2_n \leq n^{-2} x_n = y_n + 2 \frac{n + 1}{n^2} \sum_{j = 1}^{n - 1}
\frac{j^2}{(j + 1) (j + 2)} y_j, \qquad n \geq 1.
\end{equation}

We substitute~\eqref{lear}, treating the three terms separately,
into~\eqref{duncan}.  The first term in~\eqref{lear} yields the sum
\begin{eqnarray*}
\sum_{j = 1}^{n - 1} \frac{j^2}{(j + 1) (j + 2)} \frac{2 \sigma}{j^3} \sum_{k =
  1}^{j - 1} k a_k
  &=& \sum_{k = 1}^{n - 1} \sum_{k < j < n} \sigma k a_k
        \frac{2}{j (j + 1) (j + 2)} \\
  &=& \sum_{k = 1}^{n - 1} \sigma k a_k \sum_{j = k + 1}^{n - 1} \left(
        \frac{1}{j (j + 1)} - \frac{1}{(j + 1) (j + 2)} \right) \\
  &=& \sum_{k = 1}^{n - 1} \sigma k a_k \left( \frac{1}{(k + 1) (k + 2)} -
        \frac{1}{n (n + 1)} \right)
\end{eqnarray*}
and the total contribution
\begin{eqnarray}
\lefteqn{\hspace{-1.5in}\frac{2 \sigma}{n^3} \sum_{k = 1}^{n - 1} k a_k + 2
\frac{n + 1}{n^2}
  \sum_{k = 1}^{n - 1} \sigma k a_k \left( \frac{1}{(k + 1) (k + 2)} -
  \frac{1}{n (n + 1)} \right)} \nonumber \\
\label{cordelia}
  &=& 2 \sigma \frac{n + 1}{n^2} \sum_{k = 1}^{n - 1} \frac{k a_k}{(k + 1) (k +
        2)}.
\end{eqnarray}
The second term in~\eqref{lear} yields the sum
$$
\sum_{j = 1}^{n - 1} \frac{j^2}{(j + 1) (j + 2)} \frac{2 \sigma^2}{3 j^2} =
\frac{2 \sigma^2}{3} \sum_{j = 1}^{n - 1} \left( \frac{1}{j + 1} - \frac{1}{j
+ 2} \right) = \frac{2 \sigma^2}{3} \left( \frac{1}{2} - \frac{1}{n + 1}
\right)
$$
and the total contribution
\begin{equation}
\label{regan}
\frac{2 \sigma^2}{3 n^2} + 2 \frac{n + 1}{n^2} \frac{2 \sigma^2}{3} \left(
\frac{1}{2} - \frac{1}{n + 1} \right) = \frac{2 \sigma^2}{3 n^2} (1 + n + 1 -
2) = \frac{2 \sigma^2}{3 n}.
\end{equation}
Hence we find from~\eqref{duncan}
\begin{equation}
\label{macbeth}
a^2_n \leq 2 \sigma \frac{n + 1}{n^2} \sum_{k = 1}^{n - 1} \frac{k a_k}{(k +
1) (k + 2)} + \frac{2 \sigma^2}{3 n} + b^2_n + 2 \frac{n + 1}{n^2} \sum_{k =
1}^{n - 1} \frac{k^2 b^2_k}{(k + 1) (k + 2)}.
\end{equation}

We next use the following estimate of~$b_n$, whose proof we postpone.

\begin{lemma}
\label{L:A1}
For $n \geq 1$, 
$$
b_n := \| C_n(\lc n U \rc) - C(U) \|_2 \leq \left( 3 + \frac{2 \pi}{\sqrt{3}}
\right) \frac{1}{n} < \frac{6.63}{n}.
$$
\end{lemma}
\smallskip

Using this lemma in~\eqref{macbeth}, we find in analogy with~\eqref{regan}
\begin{equation}
\label{goneril}
b^2_n + 2 \frac{n + 1}{n^2} \sum_{k = 1}^{n - 1} \frac{k^2 b^2_k}{(k + 1) (k +
2)} < \frac{(6.63)^2}{n} < \frac{44}{n}
\end{equation}
and thus
\begin{equation}
\label{romeo}
a^2_n \leq 2 \sigma \frac{n + 1}{n^2} \sum_{k = 1}^{n - 1} \frac{k a_k}{(k + 1)
(k + 2)} + \left(44 + \frac{2 \sigma^2}{3} \right) \frac{1}{n}, \qquad n \geq
1.
\end{equation}

We claim that~\eqref{romeo} implies the sought estimate $a_n = O(n^{-1/2})$. 
Indeed, assume that  $n \geq 1$ and that $A > 0$ is a number such that
\begin{equation}
\label{juliet}
a_k \leq A / \sqrt{k}
\end{equation}
for $1 \leq k \leq n - 1$. 
Then, using
$k + 1 \geq [k (k + 2)]^{1 / 2}$,
\begin{eqnarray}
\sum_{k = 1}^{n - 1} \frac{k a_k}{(k + 1) (k + 2)}
  &\leq& A \sum_{k = 1}^{n - 1} \frac{k^{1 / 2}}{(k + 1) (k + 2)}
   \leq  A \sum_{k = 1}^{n - 1} \frac{1}{(k + 2)^{3 / 2}} \nonumber \\
\label{capulet}
  &\leq& A \int_0^{n - 1}\!\!\frac{dx}{(x + 2)^{3 / 2}} = 2 A \left( 2^{- 1 /
           2} - (n + 1)^{- 1 / 2} \right).
\end{eqnarray}
In particular, for $n \geq 2$,
\begin{equation}
\label{tybalt}
\frac{1}{n} \sum_{k = 1}^{n - 1} \frac{k a_k}{(k + 1) (k + 2)} \leq
\frac{1}{n} 2 A \leq 2 A (n + 1)^{- 1 / 2}
\end{equation}
and thus~\eqref{capulet} yields (trivially for $n = 1$, too)
$$
\frac{n + 1}{n} \sum_{k = 1}^{n - 1} \frac{k a_k}{(k + 1) (k + 2)} \leq 2^{1 /
2} A.
$$
Consequently, by~\eqref{romeo},
$$
n a^2_n \leq
2^{3 / 2} \sigma A + 44 + 2 \frac{\sigma^2}{3}
\leq 2^{3 / 2} \sigma A + 45.
$$
If $2^{3 / 2} \sigma A + 45 \leq A^2$, which holds for example for $A = 8$,
then this yields $n a^2_n \leq A^2$, and thus~\eqref{juliet} holds for $k =
n$, too.  By
induction, \eqref{juliet} holds for all $k \geq 1$, and we have
proved the explicit estimate 
\begin{equation}
\label{juliet8}
a_n \leq \frac{8}{\sqrt{n}}, \qquad n \geq 1.
\end{equation}

This is the desired estimate, apart from the value of the constant.  To
improve the constant, we use numerical calculations by computer.  Indeed,
for~\eqref{iago},
\begin{eqnarray*}
b^2_n &=& \sum_{i = 1}^n \int_{(i - 1) / n}^{i / n}\!\left( C(u) - C_n(i)
            \right)^2\,du \\
      &=& \int_0^1\!C(u)^2\,du - 2 \sum_{i = 1}^n C_n(i) \int_{(i - 1) /
            n}^{i / n}\!C(u)\,du + \sum_{i = 1}^n \frac{1}{n} C_n(i)^2 \\
      &=& \frac{\sigma^2}{3} - 2 \sum_{i = 1}^n C_n(i) \left(
F\left( \frac{i}{n} \right) - F\left( \frac{i - 1}{n} \right) \right) +
\frac{1}{n} \sum_{i = 1}^n C_n(i)^2,
\end{eqnarray*}
where $F(u) := u^2 \ln u - (1 - u)^2 \ln(1 - u)$ and $C_n(i)$ is
given by~\eqref{hermes}; so, given any integer~$N$, $b_n$ can be
computed exactly for $n \leq N$.  Next, for $n = 1, \ldots, N$, an
upper bound $\abar_n$ to $a_n$ can be computed recursively
from~\eqref{shylock} or, equivalently, \eqref{macbeth}, using the
already computed $\abar_k$, $k < n$, to bound $a_k$ in the right-hand
side.  (We do not know how to compute $a_n$ exactly even for $n =
3$.)  For larger~$n$, we use the estimates~\eqref{juliet} and
\refL{L:A1}.

Let
\begin{eqnarray*}
V_n     &:=& \sum_{k = 1}^n \frac{k a_k}{(k + 1) (k + 2)}, \\
\Vbar_n &:=& \sum_{k = 1}^n \frac{k \abar_k}{(k + 1) (k + 2)}, \\         
W_n     &:=& \sum_{k = 1}^n \frac{k^2 b^2_k}{(k + 1) (k + 2)}.
\end{eqnarray*}
Then for $n > N$, arguing as in~\eqref{capulet}, for any~$A$ such
that~\eqref{juliet} holds for  all~$k$,
$$
V_{n - 1} \leq \Vbar_N + \sum_{k = N + 1}^{n - 1} \frac{A}{(k + 2)^{3 / 2}}
\leq \Vbar_N + 2 A \left( (N + 2)^{- 1 / 2} - (n + 1)^{- 1 / 2} \right)
$$
and thus by~\eqref{tybalt}
\begin{equation}
\label{montague}
\sfrac{n + 1}{n} V_{n - 1} = V_{n - 1} + \sfrac{1}{n} V_{n - 1} \leq
\Vbar_N + 2 A (N + 2)^{- 1 / 2}.
\end{equation}
Similarly, with $B := \left( 3 + \frac{2 \pi}{\sqrt{3}} \right)^2 < 44$, for $n
> N$, by \refL{L:A1}, we have
\begin{eqnarray*}
\sfrac{n + 1}{n} W_{n - 1}
  &\leq& W_N + \sum_{k = N + 1}^{n - 1} \frac{B}{(k + 1) (k + 2)}
           + \frac{1}{n} \sum_{k = 1}^{n - 1} \frac{B}{(k + 1) (k + 2)} \\
  &  = & W_N + B \left( \frac{1}{N + 2} - \frac{1}{n + 1} + \frac{1}{2n} -
\frac{1}{n (n + 1)} \right) = W_N + \frac{B}{N + 2} - \frac{B}{2 n}. 
\end{eqnarray*}
Consequently, \eqref{macbeth} yields, using \refL{L:A1} again
and~\eqref{montague},
\begin{eqnarray*}
a^2_n &\leq& 2 \sigma \frac{n + 1}{n^2} V_{n - 1} + \frac{2 \sigma^2}{3 n} +
               \frac{B}{n^2} + 2 \frac{n + 1}{n^2} W_{n - 1} \\
      &\leq& \frac{1}{n} \left( 2 \sigma \Vbar_N + 4 \sigma A (N + 2)^{- 1 /
               2} + \frac{2 \sigma^2}{3} + 2 W_N + 2 B (N + 2)^{- 1} \right),
               \qquad n > N.
\end{eqnarray*}
In other words, \eqref{juliet} holds for  $k > N$, with~$A$ replaced by
\begin{equation}
\label{mercutio}
A_N := \left( 2 \sigma \Vbar_N + 4 \sigma A (N + 2)^{- 1 / 2} + \frac{2
\sigma^2}{3} + 2 W_N + 2 B (N + 2)^{- 1} \right)^{1 / 2}.
\end{equation}

For $N = 100$ we find (using {\tt Mathematica}
or {\tt Maple}), rounded to four decimal places, 
$\Vbar_{100} \doteq 1.1995$ and $W_{100} \doteq 0.3466$, and thus, taking
$A = 8$ as in~\eqref{juliet8}, $A_{100} \doteq 2.3332$.  Moreover, the
computer verifies that $n^{1 / 2} \abar_n < 1.7$ for $n \leq 100$;
thus~\eqref{juliet} holds for all $k \geq 1$ with $A = 2.34$.  Using this
value in~\eqref{mercutio} we find $A_{100} \doteq 1.9976$, and the
theorem is proved.
\end{proof}

\begin{remark}
The sequence $n^{1 / 2} \abar_n$ seems to increase slowly.  For $n = 100$ the
value is (rounded to four decimal places) $1.6018$, and hence the bound in
\refT{T:d2} cannot be much improved using the present method based
on~\eqref{shylock}.
\end{remark}

It remains to prove \refL{L:A1} above.

\begin{proof}[Proof of \refL{L:A1}]
Let $I_i := \{u: \lc n u \rc = i\} = ((i - 1) / n, i / n]$.  Thus $I_1,
\ldots, I_n$ form a partition of $(0, 1]$.  We choose a point  $t_i \in
\bar{I_i}$ for each~$i$ (where the bar here indicates closure) and define
$$
\Ct_n(u) := C(t_{\lc n u \rc}),
$$
i.e.,\ $\Ct(u) = C(t_i)$ when $u \in I_i$.  By Minkowski's inequality, 
\begin{equation}
\label{hamlet}
b_n \leq \| C_n(\lc n U \rc) - \Ct_n(U) \|_2 + \| \Ct_n(U) - C(U) \|_2.
\end{equation}

To estimate the second term in~\eqref{hamlet}, note that for $u \in I_i$,
$$
|\Ct_n(u) - C(u)| = |C(t_i) - C(u)| \leq \int_{I_i}\!|C'(x)|\,dx.
$$
The \CS\ inequality yields
$$ 
|\Ct_n(u) - C(u)|^2 \leq \frac{1}{n} \int_{I_i}\!|C'(x)|^2\,dx, \qquad
u \in I_i,
$$
and thus (for any choice of  $t_i \in \bar{I_i}$),
\begin{eqnarray}
\| \Ct_n(U) - C(U) \|^2_2
  &  = & \sum_{i = 1}^n \int_{I_i}\!|\Ct_n(u) - C(u)|^2\,du
           \nonumber \\
  &\leq& \sum_{i = 1}^n \frac{1}{n^2} \int_{I_i}\!|C'(x)|^2\,dx \nonumber \\
\label{ophelia}
  &  = & \frac{1}{n^2} \int_0^1 |C'(x)|^2\,dx.
\end{eqnarray}
We have
$$
C'(x) = 2 \ln x - 2 \ln(1 - x)
$$
and find
$$
\int_0^1\!\left[ \ln (1 - x) \right]^2\,dx = \int_0^1\!(\ln x)^2\,dx =
\int_0^{\infty}\!y^2 e^{-y}\,dy = 2
$$
and
\begin{eqnarray*}
\int_0^1\![\ln x] [\ln(1 - x)]\,dx
  &=& \sum_{k = 1}^{\infty} \frac{1}{k} \int_0^1\!x^k\,|\ln x|\,dx
   =  \sum_{k = 1}^{\infty} \frac{1}{k} \int_0^{\infty}\!e^{- k y} y e^{-
        y}\,dy \\
  &=& \sum_{k = 1}^{\infty} \frac{1}{k} \frac{1}{(k + 1)^2}
   =  \sum_{k = 1}^{\infty} \left( \frac{1}{k (k + 1)} - \frac{1}{(k + 1)^2}
        \right) \\
  &=& 1 - \left( \frac{\pi^2}{6} - 1 \right) = 2 - \frac{\pi^2}{6};
\end{eqnarray*}
consequently,
$$
\int_0^1\!|C'(x)|^2\,dx = 8 \int_0^1\!(\ln x)^2\,dx - 8 \int_0^1\![\ln x]
[\ln(1 - x)]\,dx = \frac{4 \pi^2}{3}.
$$
Hence~\eqref{ophelia} yields
\begin{equation}
\label{gertrude}
\| \Ct_n(U) - C(U) \|_2 \leq \frac{1}{n} \| C'(U) \|_2 \leq \frac{2
\pi}{\sqrt{3} n}.
\end{equation}

For the first term in~\eqref{hamlet}, let  us first assume that $n \geq 2$.
For $u \in I_i$ we have
\begin{eqnarray*}
C_n(\lc n u \rc) - \Ct_n(u)
  &=& C_n(i) - C(t_i) \\
  &=& - \sfrac{1}{n} + \sfrac{1}{n} (\mu_{i - 1} + \mu_{n - i} - \mu_n) - 2 t_i
        \ln t_i - 2 (1 - t_i) \ln (1 - t_i).
\end{eqnarray*}
For $i \leq \lc n / 2 \rc$ we choose $t_i = i / n$.  This yields,
using~\eqref{clymenestra} and~\eqref{helen},
\begin{eqnarray}
C_n(i) - C(t_i)
  &\leq& \sfrac{1}{n} \left[ -1 + 2 i \ln i + (2 \gamma - 4) i + 3
           \right. \nonumber \\
  &    & \qquad + 2 (n - i + 1) \ln(n - i) + (2 \gamma - 4) (n - i) + 2
           \gamma + 1 + \sfrac{1}{n - i} \nonumber \\
  &    & \qquad \left. - 2 (n + 1) \ln n - (2 \gamma - 4) n - 2 \gamma - 2 i
           \ln(\sfrac{i}{n}) - 2 (n - i) \ln(\sfrac{n - i}{n}) \right]
           \nonumber \\
  &  = & \sfrac{1}{n} \left[ 2 \ln(\sfrac{n - i}{n}) + 3 + \sfrac{1}{n - i}
           \right]
   \leq  \sfrac{1}{n} \left[ 3 - \sfrac{2 i}{n} + \sfrac{1}{n - i} \right]
           \nonumber \\
\label{yorrick}
  &\leq& \sfrac{3}{n}.
\end{eqnarray}
In the opposite direction, by~\eqref{clymenestra} and~\eqref{helen}, still for
$i \leq \lc n / 2 \rc$,
\begin{eqnarray*}
C_n(i) - C(t_i)
  &\geq& \sfrac{1}{n} \left[ -1 + 2 i \ln i + (2 \gamma - 4) i + 2 \right. \\
  &    & \qquad + 2 (n - i + 1) \ln(n - i) + (2 \gamma - 4) (n - i) + 2 \gamma
           \\
  &    & \qquad \left. - 2 (n + 1) \ln n - (2 \gamma - 4) n - 2 \gamma - 1 -
           \sfrac{1}{n} \right.\\
  &    & \qquad \left. - 2 i \ln(\sfrac{i}{n}) - 2 (n - i) \ln(\sfrac{n -
           i}{n}) \right] \\
  &  = & \sfrac{1}{n} \left[ 2 \ln(\sfrac{n - i}{n}) - \sfrac{1}{n}
           \right]
   \geq  \sfrac{1}{n} \left[ 2 \ln(\sfrac{1}{3}) - \sfrac{1}{2} \right]
   \geq  - \sfrac{3}{n}.
\end{eqnarray*}
Consequently, for $i \leq \lc n / 2 \rc$,
\begin{equation}
\label{polonius}
|C_n(i) - C(t_i)| \leq 3 / n.
\end{equation}
For $i > \lc n / 2 \rc$ we choose instead $t_i = (i - 1) / n = 1 - t_{n + 1 -
i}$.  The symmetries of~$C_n$ and~$C$ then yield $C_n(i) - C(t_i) = C_n(n + 1
- i) - C(t_{n + 1 - i})$, and since $n + 1 - i \leq n / 2$, \eqref{polonius}
shows that $|C_n(i) - C(t_i)| \leq 3 / n$ for $i > \lc n / 2 \rc$, too, i.e.,\
\eqref{polonius} holds for all $i \leq n$.  In other words,
$$
|C_n(\lc n u \rc) - \Ct_n(u)| = |C_n(\lc n u \rc) - C(t_{\lc n u \rc})| \leq 3
/ n
$$
for all $u \in (0, 1]$; in particular, $\| C_n(\lc n U \rc) -
\Ct_n(U) \|_2 \leq 3 / n$  for all $n \geq 2$.  This holds trivially for
$n = 1$, too, for any choice of $t_1$, and
together with~\eqref{hamlet} and~\eqref{gertrude} yields the result.
\end{proof}

\begin{remark}
Define 
$$
c^* := \sup \{n^{1 / 2} d_2(Y_n, Y):\ n \geq 1\},
$$
so that, by \refT{T:d2}, $c^* \leq 2$.  Conversely,
$$
c^* \geq 2^{1 / 2} d_2(Y_2, Y) = 2^{1 / 2} \| Y \|_2 = \sigma \sqrt{2} >
0.9168;
$$
thus the constant~$2$ in \refT{T:d2} is no more than about twice the optimal
value.

Although we do not know the exact value of
$d_2(Y_n, Y)$ for any
$n > 2$, one can in principle for any~$n$ and~$m$ compute the exact
distributions of~$Y_n$ and~$Y_m$ and thus $d_2(Y_n, Y_m)$.  We have done this
for some $m, n \leq 50$ using {\tt Mathematica} and {\tt Maple}.
The results are consistent with a decay of the type
$d_2(Y_n, Y) \sim c n^{- 1 / 2}$ with $c \approx 1$, but our data are too few
to be conclusive.
\end{remark}

\section{Bounding $d_p(Y_n, Y)$} \label{S:dp}

In this section we extend \refT{T:d2} and show that $d_p(Y_n, Y) = O(n^{- 1 /
2})$ for every~$p$.  In contrast to the style of \refS{S:d2}, we will make no
attempt to keep constants small, nor to keep track of them explicitly.

\begin{theorem}\label{T:dp}
For every $p \geq 1$, there exists a constant $c_p < \infty$ such that
$$
d_p(Y_n, Y) \leq c_p / \sqrt{n}, \qquad n \geq 1.
$$
\end{theorem}

\begin{proof}
Since $d_p \leq d_q$ when $p \leq q$, it suffices to consider integer $p \geq
2$.  The case $p = 2$ is \refT{T:d2} (with $c_2 = 2$), so we assume further
that $p \geq 3$.  We use induction on~$p$ and assume that the result holds for
smaller positive integer values of~$p$.

Let $Y, Y_n, Y^*, Y^*_n, U$ be as in \refS{S:d2}, and note that for every
$p \geq 1$,
\begin{equation}
\label{p1}
\| Y_n - Y \|_p = \| Y^*_n - Y^* \|_p = d_p(Y_n, Y), \qquad n \geq 0,
\end{equation}
by the fact~\cite{Cambanis} that there is an optimal coupling for~$d_2$
that is optimal for every~$d_p$.  Using the notation of
\refS{S:d2}, we have, for $n \geq 1$,
\begin{equation}
\label{p2}
d_p(Y_n, Y) \leq \| \Yt_n - \Yt \|_p = \| W_1 + W_2 + W_3 \|_p.
\end{equation}
We use a simple lemma to estimate this.

\begin{lemma}
\label{L:three}
Let $Z_1$, $Z_2$, and $Z_3$ be three independent random variables, and let $p
\geq 2$ be an integer.  Then 
$$
\EE\,|Z_1 + Z_2 + Z_3|^p \leq \EE\,|Z_1|^p + \EE\,|Z_2|^p + \left( \|Z_1\|_{p
- 1} + \|Z_2\|_{p - 1} + \|Z_3\|_p \right)^p.
$$
\end{lemma}

\begin{proof}
By the binomial theorem and independence,
$$
\EE\,|Z_1 + Z_2 + Z_3|^p \leq \EE \left( |Z_1| + |Z_2| + Z_3| \right)^p =
\sum_{j, k, l} \binom{p}{j, k, l} \left( \EE\,|Z_1|^j \right)
\left( \EE\,|Z_2|^k \right) \left( \EE\,|Z_3|^l \right).
$$
If $j \leq p - 1$ and $k \leq p - 1$ we estimate
$\EE\,|Z_1|^j = \|Z_1\|^j_j \leq \|Z_1\|^j_{p - 1}$ (which holds also for $j =
0$, disregarding the central expression) and similarly $\EE\,|Z_2|^k \leq
\|Z_2\|^k_{p - 1}$ and $\EE\,|Z_3|^l \leq \|Z_3\|^l_p$.  Hence all terms in
the sum, except $\EE\,|Z_1|^p$ and $\EE\,|Z_2|^p$, are bounded by the
corresponding terms in the trinomial expansion of $\left( \|Z_1\|_{p - 1} +
\|Z_2\|_{p - 1} + \|Z_3\|_p \right)^p$.
\end{proof}

Conditional on $U = u$, the three variables
$W_1$, $W_2$, and $W_3$ are independent, so the lemma is applicable.  Fix $u
\in (0, 1)$ and let $i = \lc n u \rc$, so $1 \leq i \leq n$.  Then, given
$U = u$, $W_1 = \frac{i - 1}{n} Y_{i - 1} - u Y$ and thus, for any $q \geq
1$,
\begin{eqnarray}
\EE \left( |W_1|^q\,|\,U = u \right)^{1 / q}
  &  = & \| \sfrac{i - 1}{n} Y_{i - 1} - u Y \|_q \nonumber \\
  &\leq& \| \sfrac{i - 1}{n} (Y_{i - 1} - Y) \|_q + |\sfrac{i - 1}{n} - u|\,
           \|Y\|_q \nonumber \\
\label{w1}
  &\leq& \sfrac{i - 1}{n} d_q(Y_{i - 1}, Y) + \sfrac{1}{n} \|Y\|_q.
\end{eqnarray}
Similarly, 
\begin{equation}
\label{w2}
\EE \left( |W_2|^q\,|\,U = u \right)^{1 / q} \leq \sfrac{n - i}{n}
d_q(Y_{n - i}, Y) + \sfrac{1}{n} \|Y\|_q.
\end{equation}
Further, given $U = u$, $W_3 = C_n(i) - C(u)$ is a constant, for which we
use the simple estimate (from Proposition~3.2 in~\cite{Roesler})
\begin{equation}
\label{w3}
|W_3| = |C_n(\lc n u \rc) - C(u)| \leq \sfrac{6}{n} \ln n + O \left( n^{-
1} \right) = O \left( n^{- 1 / 2} \right).
\end{equation}

We first use~\eqref{w1} with $q = p - 1$ together with the induction
hypothesis $d_{p - 1}(Y_{i - 1}, Y) \leq  c_{p - 1} (i - 1)^{- 1 / 2}$, $i
\geq 2$, to obtain (also for $i = 1$)
$$
\EE \left( |W_1|^{p - 1}\,|\,U = u \right)^{1 / (p - 1)} \leq c_{p - 1}
\frac{(i - 1)^{1 / 2}}{n} + \frac{1}{n} \|Y\|_{p - 1} \leq b_1 n^{- 1 / 2},
$$
where $b_1$, like $b_2, b_3, b_4$
below, denotes some constant depending
on~$p$ only.  By similar argument using~\eqref{w2} and~\eqref{w3}, we obtain
$$
\EE \left( |W_1|^{p - 1}\,|\,U = u \right)^{1 / (p - 1)} + \EE \left(
|W_2|^{p - 1}\,|\,U = u \right)^{1 / (p - 1)} + \EE \left( |W_3|^p\,|\,U
= u \right)^{1 / p} \leq b_2 n^{- 1 / 2}.
$$
Hence, using~\eqref{w1} and~\eqref{w2} for $q = p$, too, \refL{L:three} yields
\begin{eqnarray*}
\EE \left( |W_1 + W_2 + W_3|^p\,|\,U = u \right)
  &\leq& \left( \frac{i - 1}{n} d_p(Y_{i - 1}, Y) + b_3 \frac{1}{n}
           \right)^p \\
  &    & \qquad + \left( \frac{n - i}{n} d_p(Y_{n - i}, Y) + b_3 \frac{1}{n}
           \right)^p + b^p_2 n^{- p / 2}.
\end{eqnarray*}
Taking the average over all  $u \in (0, 1)$ we finally find the recursive
estimate
\begin{eqnarray}
d_p(Y_n, Y)^p
  &\leq& \EE\,|W_1 + W_2 + W_3|^p = \EE\,\EE \left( |W_1 + W_2 +
           W_3|^p\,|\,U \right) \nonumber \\
\label{dpind}
  &\leq& \frac{2}{n} \sum_{j = 0}^{n - 1} \left( \frac{j}{n} d_p(Y_j, Y) +
           b_3 \frac{1}{n} \right)^p + b^p_2 n^{- p / 2}.
\end{eqnarray}

The proof is now completed by another induction, this one on~$n$.  Suppose
that $d_p(Y_j, Y) \leq c j^{- 1 / 2}$ for $1 \leq j \leq n - 1$. 
Then~\eqref{dpind} yields
\begin{eqnarray}
d_p(Y_n, Y)^p
  &\leq& \frac{2}{n} \sum_{j = 1}^{n - 1} \left( c j^{1 / 2} n^{- 1} +
           b_3 n^{- 1} \right)^p + \frac{2}{n} b^p_3 n^{- p} + b^p_2
           n^{- p / 2} \nonumber \\
  &\leq& \frac{2}{n} (c + b_3)^p \sum_{j = 1}^{n - 1} j^{p / 2} n^{- p} +
           b_4 n^{- p / 2} \nonumber \\
  &\leq& 2 (c + b_3)^p \int_0^1\!x^{p / 2} n^{- p / 2}\,dx + b_4 n^{- p /
           2} \nonumber \\
\label{dp2}
  &  = &  \left[ 2 (c + b_3)^p \frac{1}{(p / 2) + 1} + b_4 \right] n^{- p
           / 2}.
\end{eqnarray}
Since $p \geq 3$, we have $\frac{2}{(p / 2) + 1} = \frac{4}{p + 2} < 1$, and
thus, if~$c$ is sufficiently large,
$$
\frac{4}{p + 2} (c + b_3)^p + b_4 \leq c^p.
$$
For such~$c$, \eqref{dp2} yields $d_p(Y_n, Y)^p \leq \left( c n^{- 1 / 2}
\right)^p$, which completes both inductions and the proof.
\end{proof}

Note that the arguments used above for $p \geq 3$ do not work for $p = 2$, so
we need both the proof here and the proof in \refS{S:d2}.

\section{Lower bounds for $d_p(Y_n, Y)$} \label{S:lower}

We
do not know whether the upper bounds $O(n^{- 1 / 2})$ proved in the
preceding two sections are sharp.  We give in this section two simple lower
bounds.

First, $d_p(Y_n, Y) = \Omega(n^{- 1})$ for every~$p$ by the following general
result.

\begin{proposition}
Let $W, W_1, W_2, \ldots$ be random variables such that~$W$ has an absolutely
continuous distribution while, for each $n \geq 1$ and some constant~$b_n$,
$n (W_n - b_n)$ is integer-valued.  Then, for each $1 \leq p < \infty$,
$d_p(W_n, W) =
\Omega(1 / n)$.  More precisely,
\begin{equation}
\label{rome}
\liminf_{n \to \infty} n\,d_p(W_n, W) \geq \sfrac{1}{2} (p + 1)^{- 1 / p}.
\end{equation}
\end{proposition}

\begin{proof}
Let 
$V_n := \{n (W - b_n)\}$, where $\{x\} := x - \lf x \rf$ denotes the
fractional part of~$x$.  For any coupling of~$W$ and~$W_n$,
$$
|W - W_n| = \sfrac{1}{n} |n (W - b_n) - n (W_n - b_n)| \geq \sfrac{1}{n}
h(V_n),
$$
where $h(x) := \min(x, 1 - x)$, $0 \leq x \leq 1$, and thus
\begin{equation}
\label{venice}
d_p(W, W_n) \geq \sfrac{1}{n} \|h(V_n)\|_p.
\end{equation}
We regard~$V_n$ as a random variable taking values in $\RR/\ZZ \cong
\TT$, and find that its distribution, $\nu_n$
say, has Fourier coefficients
$$
\nuh_n(k) = \EE \left( e^{- 2 \pi i k V_n} \right) = e^{2 \pi i k n b_n} \phi(-
2 \pi k n),
$$
where~$\phi$ is the characteristic function of~$W$.  In particular,
$|\nuh_n(k)| = |\phi(2 \pi k n)|$.  By our
hypothesis on~$W$ and the Riemann--Lebesgue lemma, $\phi(x) \to 0$ as $x
\to \pm \infty$.  Thus, for each fixed $k \neq 0$, $\nuh_n(k) \to 0$ as $n \to
\infty$.  This implies that~$\nu_n$ converges weakly (as measures on~$\TT$) to
the uniform distribution, i.e.,\ $V_n\,\Lto\,U$ where $U \sim \mbox{unif}(0,
1)$.  Consequently, as $n \to \infty$,
$$
\|h(V_n)\|^p_p = \EE\,h(V_n)^p \to \EE\,h(U)^p = 2 \int_0^{1 / 2}\!x^p\,dx =
2^{- p} / (p + 1),
$$
which together with~\eqref{venice} leads to~\eqref{rome}.  The proof of the
proposition is completed by observing $d_p(W_n, W) > 0$ for every~$n$, because
$W_n \Lnoteq W$.
\end{proof}

Note that, in contrast to the asymptotic result~\eqref{rome}, there is no
positive lower bound to $d_p(W_n, W)$ for a fixed~$n$ without further
assumptions.  Hence the implicit constant in $\Omega(1 / n)$ in the theorem
depends on the variables $W, W_1, \ldots$.

For $p \geq 2$ we can improve this lower bound by a logarithmic factor by
using the known variance of~$Y_n$.

\begin{theorem}
If $2 \leq p < \infty$, then
$$
d_p(Y_n, Y) \geq d_2(Y_n, Y) = \Omega(\sfrac{\ln n}{n}).
$$
\end{theorem}

\begin{proof}
Recall that~$Y$ and~$Y_n$ have mean~$0$ and that $\Var\,Y = \sigma^2$ while
by~\eqref{poseidon}
$$
\Var\,Y_n = \sigma^2 - 2 \sfrac{\ln n}{n} + O(n^{- 1})
$$
and thus
$$
\| Y_n \|_2 = (\Var\,Y_n)^{1 / 2} = \sigma - \sfrac{1}{\sigma} \sfrac{\ln
n}{n} + O(n^{- 1}).
$$
Consequently, for the $d_2$-optimal coupling of~$Y$ and~$Y_n$, by Minkowski's
inequality,
$$
d_2(Y_n, Y) = \| Y_n - Y \|_2 \geq \| Y \|_2 - \| Y_n \|_2 = \sigma^{- 1}
\sfrac{\ln n}{n} + O(n^{- 1}).
$$
\end{proof}

We still have a gap between $(\ln n) / n$ and $n^{- 1 / 2}$.

\begin{remark}
It
can be shown that $\EE\,Y^m_n = \EE\,Y^m + O \left( \frac{\ln n}{n}
\right)$, $n \geq 2$, holds also for $m = 3, 4, \ldots$; cf.\ the formulas for
moments and cumulants by Hennequin~\cite{H91}. Hence we do not get better lower
bounds for~$d_p$ by considering higher moments.
\end{remark}

\section{The \KSm\ distance} \label{S:KS}

Recall that the \KSm\ distance $\dks(F, G)$ between two distributions is
defined as $\sup_{x \in \RR} |P(X \leq x) - P(Y \leq x)|$, when $X \sim F$ and
$Y \sim G$.  We will in this case also write $\dks(X, Y)$.

To obtain upper bounds for $\dks(Y_n, Y)$, we combine the bounds above for
$d_p(Y_n, Y)$ with the following simple general result and the
fact~\cite{SJFill1} that~$Y$ has a bounded density function.

\begin{lemma}
\label{L:KS1}
Suppose that~$X$ and~$Y$ are two random variables such that~$Y$ is absolutely
continuous with a bounded density function~$f$.  If $M := \sup_{y \in \RR}
|f(y)|$ and $1 \leq p < \infty$, then
$$
\dks(X, Y) \leq (p + 1)^{1 / (p + 1)} \left( M d_p(X, Y) \right)^{p / (p +
1)}.
$$
\end{lemma}

\begin{proof}
Consider an optimal $d_p$-coupling of~$X$ and~$Y$.  Then, for $x \in \RR$ and
$\eps > 0$, denoting the distribution functions of~$X$ and~$Y$ by~$F_X$
and~$F_Y$,
\begin{eqnarray*}
F_X(x) = P(X \leq x)
  &\leq& P(Y \leq x) + P(x < Y \leq x + \eps) + P(Y - X > \eps) \\
  &\leq& F_Y(x) + M \eps + P(Y - X > \eps).
\end{eqnarray*}
Similarly,
\begin{eqnarray*}
F_Y(x)
  &\leq& P(X \leq x) + P(x - \eps < Y \leq x) + P(X - Y > \eps) \\
  &\leq& F_X(x) + M \eps + P(X - Y > \eps).
\end{eqnarray*}
Consequently,
$$
\Delta(x) := |F_X(x) - F_Y(x)| \leq M \eps + P(|X - Y| > \eps)
$$
and thus
\begin{align*}
d_p(X, Y)^p
  &  =  \EE\,|X - Y|^p = \int_0^{\infty}\!p \eps^{p - 1} P(|X - Y| >
           \eps)\,d\eps \\
  &\geq \int_0^{\Delta(x) / M}\!p \eps^{p - 1} (\Delta(x) - M \eps)\,d\eps
           \\
  &  =  \sfrac{1}{p + 1} \Delta(x)^{p + 1} M^{- p}. \qedtag
\end{align*}
\noqed
\end{proof}

\begin{theorem}
\label{T:KS1}
For every $\eps > 0$,
$$
\dks(Y_n, Y) = O \left( n^{\eps - (1 / 2)} \right).
$$
\end{theorem}

\begin{proof}
By~\cite{SJFill1}, $Y$ has a bounded density function, so \refL{L:KS1} and
\refT{T:dp} yield, for every fixed $1 \leq p < \infty$,
$$
\dks(Y_n, Y) = O \left( d_p(Y_n, Y)^{p / (p + 1)} \right) = O \left( n^{- p
/ [2 (p + 1)]} \right).
$$
The result follows by choosing~$p$ so large that $\frac{p}{2 (p + 1)} >
\frac{1}{2} - \eps$.
\end{proof}

To  get an explicit bound we take $p = 2$ in \refL{L:KS1} and use
\refT{T:d2}.  This yields the bound $3^{1 / 3} \left(2 M n^{- 1 / 2} \right)^{2
/ 3}$, and we know $M < 16$ from Theorem~3.3 of~\cite{SJFill1}.
Hence,

\begin{theorem}
\label{T:KS2}
For $n \geq 1$,
$$
\dks(Y_n, Y) \leq (12 M^2)^{1 / 3} n^{- 1 / 3} < (3072 / n)^{1 / 3} < 15\,
n^{- 1 / 3}. \nopf
$$
\end{theorem}
\noindent
Numerical evidence~\cite{TanH} suggests that $M < 1$, which would
give a bound $2.3\,n^{- 1 / 3}$.

We conjecture that \refT{T:KS1} holds with $\eps = 0$, too, i.e.,\ that
$\dks(Y_n, Y) = O(n^{- 1 / 2})$.  Even if this were proved, it is not clear
what the right order of decay is; the best lower bound we can prove is
$\Omega(n^{- 1})$.
\begin{theorem}
\label{T:KSlower}
$$
\dks(Y_n, Y) > \frac{1}{8 (n + 1)}, \qquad n \geq 1.
$$
\end{theorem}

Again, the lower bound follows from quite general considerations.  In this
case we use the following lemma.

\begin{lemma}
\label{L:KS2}
Suppose that~$Y$ and~$Z$ are two random variables such that~$Y$ has a
continuous distribution while $a (Z - b)$ is integer-valued for some real
numbers $a > 0$ and~$b$.  If $\sigma^2_Z := \Var\,Z < \infty$, then
$$
\dks(Y, Z) \geq 1 / (12 a \sigma_Z + 8).
$$
\end{lemma}

\begin{proof}
For any $x \in \RR$ and $\delta > 0$,
$$
F_Z(x + \delta) - F_Y(x + \delta) + F_Y(x - \delta) - F_Z(x - \delta) \leq 2
\dks(Y, Z).
$$
Letting $\delta \to 0$ we find, since~$Y$ is continuous,
$$
P(Z = x) \leq 2 \dks(Y, Z).
$$
The result now follows from the following lemma applied to $a (Z - b)$.
\end{proof}

\begin{lemma}
If~$Z$ is an integer-valued random variable with finite variance
$\sigma^2_Z$, then
$$
\sup_n P(Z = n) \geq 1 / (6 \sigma_Z + 4).
$$
\end{lemma}

\begin{proof}
Let $\mu := \EE\,Z$ and $m := \lc \sfrac{3}{2} \sigma_Z \rc$.  By Chebyshev's
inequality,
$$
P(|Z - \mu| \geq m) \leq \frac{\sigma^2_Z}{m^2} \leq \frac{4}{9} < \frac{1}{2}
$$
and thus
$$
P(\mu - m < Z < \mu + m) > 1 / 2.
$$
The interval $(\mu - m, \mu + m)$ contains at most $2 m$ integers, and thus it
must contain an integer~$n$ such that
\begin{equation*}
P(Z = n) \geq \frac{1}{2 m} P(\mu - m < Z < \mu + m) > \frac{1}{4 m} >
\frac{1}{6 \sigma_Z + 4}.
\qedtag
\end{equation*}
\noqed
\end{proof}

\begin{proof}[Proof of \refT{T:KSlower}]
We apply \refL{L:KS2} with $a = n$ and observe that
\begin{equation}
\label{ravenna}
\sigma_{Y_n} := (\Var\,Y_n)^{1 / 2} < \sigma = (\Var\,Y)^{1 / 2}
\end{equation}
and that $12 \sigma \doteq 7.8 < 8$.  Indeed, \eqref{ravenna} is trivial for
$n = \mbox{$1$ or~$2$}$ and easily verified for $3 \leq n \leq 6$, while
for $n \geq 7$ it holds because then, by~\eqref{poseidon}
and~\eqref{artemis},
\begin{align*}
\sigma^2 - \Var\,Y_n
  &= - 4 \frac{\pi^2}{6} + 4 \left( 1 + \frac{1}{n} \right)^2 H^{(2)}_n + 2
        \frac{n + 1}{n^2} H_n - \frac{13}{n} \\
  &> - 4 \sum_{k = n + 1}^{\infty} k^{- 2} + \frac{8}{n} H^{(2)}_n +
        \frac{2}{n} H_n - \frac{13}{n} \\
  &> \sfrac{1}{n} \left( -4 + 8 H^{(2)}_n + 2 H_n - 13 \right) > 0.
\qedtag
\end{align*}   
\noqed
\end{proof}

\section{Approximating the density of~$Y$}
\label{S:local}

It was shown in~\cite{SJFill1} that the density~$f$ of~$Y$ is infinitely
differentiable, with all derivatives rapidly decaying.  In particular, the
derivative~$f'$ is bounded; Theorem~3.3 of~\cite{SJFill1} gives the explicit
bound
$$
M' := \sup_{x \in \RR} |f'(x)| < 2466.
$$
(This is not very sharp; the true value seems to be less than~$2$.)
The bounds above on the \KSm\ distance then imply the following local result.

\begin{theorem}
\label{T:f}
For any $x \in \RR$ and $\delta > 0$,
$$
\left| \frac{F_n(x + \frac{\delta}{2}) - F_n(x - \frac{\delta}{2})}{\delta}
- f(x) \right| \leq \frac{(96 M^2)^{1 / 3}}{\delta n^{1 / 3}} + \frac{M'}{4}
\delta.
$$
In particular, for any $\Mbar \geq M$ and $\Mbar' \geq M'$, choosing 
$\delta = \delta_n := 2 \left( 96 \Mbar^2 (\Mbar')^{- 3} \right)^{1 / 6} n^{- 1
/ 6}$ yields
\begin{equation}
\label{pisa}
\left| \frac{F_n(x + \frac{\delta_n}{2}) - F_n(x -
\frac{\delta_n}{2})}{\delta_n} - f(x) \right| \leq \left( 96 \Mbar^2
(\Mbar')^3 \right)^{1 / 6} n^{- 1 / 6}.  
\end{equation}
\end{theorem}

The choices $\Mbar = 16$ and $\Mbar' = 2466$ provided by~\cite{SJFill1} yield
the bound  $268\,n^{- 1 / 6}$ in~\eqref{pisa}.  If $\Mbar = 1$ and
$\Mbar' = 2$ could be proven to be legitimate, we could reduce the bound to
$3.03\,n^{- 1 / 6}$.

\begin{proof}
By \refT{T:KS2},
$$
\left| F_n \left( x + \sfrac{\delta}{2} \right) - F_n \left( x -
\sfrac{\delta}{2} \right) - \left( F \left (x + \sfrac{\delta}{2} \right) - F
\left( x - \sfrac{\delta}{2} \right) \right) \right| \leq 2 \dks(Y_n, Y) \leq 2
(12 M^2)^{1 / 3} n^{- 1 / 3},
$$
while
\begin{eqnarray*}
\left| F \left( x + \frac{\delta}{2} \right) - F \left( x - \frac{\delta}{2}
  \right) - \delta f(x) \right|
  &  = & \left| \int_{- \delta / 2}^{\delta / 2}\!(f(x + y) - f(x))\,dy
           \right| \\
  &\leq& \int_{- \delta / 2}^{\delta / 2}\!M' |y|\,dy = \frac{M'}{4} \delta^2.
\end{eqnarray*}
The first estimate follows, and the second is an immediate consequence.
\end{proof}

\refT{T:f} yields a simple method to numerically calculate the unknown
density~$f$ up to any given accuracy.  For an application, see~\cite{DFN}. 
(In~\cite{DFN}, a preliminary version of \refT{T:f} with larger constants is
used.)  Note, however, that the convergence is slow and that it seems
impractical to obtain high precision by this method.  Other, potentially
more powerful, methods to calculate~$f$ numerically are discussed
in~\cite{SJFill2}.

\begin{problem}
Does a local limit theorem hold in the form that
$$
\left| n P(X_n = k) - f \left( \frac{k - \mu_n}{n} \right) \right| = \left| n P
\left( Y_n = \frac{k - \mu_n}{n} \right) - f \left( \frac{k - \mu_n}{n}
\right) \right| \to 0,
$$
perhaps uniformly in $k \in \ZZ$, as $n \to \infty$?
\end{problem}

\section{Bounds on moment generating functions}
\label{S:mgf}

\Roesler~\cite{Roesler} proved that the moment generating functions
$\EE\,e^{\gl Y_n}$ are bounded for fixed~$\gl$, and thus
$\EE\,e^{\gl Y_n} \to \EE\,e^{\gl Y}$ as $n \to \infty$. 
\Roesler\ did not make his estimates explicit, but 
his method can be used to obtain explicit bounds.
For the limit variable~$Y$, this was done in~\cite{SJFill2},
where we obtained by \Roesler's method (with some refinements) the following
explicit estimates for the moment generating function of~$Y$:
Let $L_0 \doteq 5.018$ be the largest root of $e^L = 6 L^2$; then 
\begin{equation}
\label{homer}
\psi_Y(\gl) := \EE\,e^{\gl Y} \leq 
  \begin{cases}
    e^{1.25 \gl^2},                    & \gl \leq -0.62, \\
    e^{0.5 \gl^2}, & -0.62 \leq \gl \leq 0,\\
 	e^{\gl^2}, & 0 \leq \gl \leq 0.42, \\
    e^{12 \gl^2},                  & 0.42 \leq \gl \leq L_0, \\
    e^{2 e^{\gl}},                 & L_0 \leq \gl.
   \end{cases}
\end{equation}
In particular, $\EE\,e^{\gl Y} \leq \exp \left( \max \left( 12 \gl^2, 2
e^{\gl} \right) \right)$ for all $\gl \in \RR$.

The constants in~\eqref{homer} are
not sharp, but the doubly exponential growth as $\gl \to +\infty$ is
correct:\ it was also shown in~\cite{SJFill2} that $\psi_Y(\gl) \geq
\exp \left( \gamma \gl^{-1} e^{\gl} \right)$ for all large~$\gl$
whenever $\gamma < 2 / e$.
  
In this section we will establish
similar bounds for $\EE\,e^{\gl Y_n}$.
For simplicity we first consider the slight shrinkage
$$
\Yh_n := \frac{n}{n + 1} Y_n = \frac{X_n - \mu_n}{n + 1}
$$
of $Y_n$; in  particular, $\Yh_0:= X_0-\mu_0=0$.
We then have the following simple result. 
\begin{theorem}
\label{T:simple}
$\EE\,e^{\gl \Yh_n} \uparrow \EE\,e^{\gl Y}$ as $n\uparrow \infty$.
Hence, for any $n\ge0$, 
$\EE\,e^{\gl \Yh_n}\le \EE\,e^{\gl Y}$,
and in particular the upper bounds on $\EE\,e^{\gl Y}$
in~\eqref{homer} above apply also to
$\EE\,e^{\gl \Yh_n}$.
\end{theorem}

\begin{proof}
It is well known that the number $X_n$ of \Quicksort{} comparisons
has the same distribution as the internal path length of a random 
binary search tree (under the random permutation model) with $n$ internal
nodes---see, e.g.,~\cite[Section 6.2.2]{Knuth3}. Moreover, it was shown by 
R\'{e}gnier~\cite{Reg} that when 
$X_n$ is reinterpreted as the internal path length of an evolving
random binary search tree after~$n$ keys have been inserted, the process 
$(\Yh_n)_{n \geq 0}$ is a martingale, which is $L^2$-bounded and thus
converges a.s.\ and in $L^2$ to some limit $Y$. It follows that also
$Y_n\to Y$ a.s., and thus in distribution; hence this random variable $Y$
is (a realization of) the same $Y$ as above.

The martingale property can be written 
$\Yh_n = \EE(\Yh_{n+1}|\Fc_n)$, for the appropriate $\sigma$-field
$\Fc_n$. Since $x\mapsto e^{\gl x}$ is convex, it now follows by
Jensen's inequality for conditional expectations that 
$e^{\gl \Yh_n} \le \EE(e^{\gl\Yh_{n+1}}|\Fc_n)$; 
and thus, taking expectations,
$\EE e^{\gl \Yh_n} \le \EE e^{\gl\Yh_{n+1}}$. 

By the same argument, $\EE e^{\gl \Yh_n} \le \EE e^{\gl Y}$ for each $n\ge0$,
which together with Fatou's lemma yields
$\EE e^{\gl \Yh_n} \to \EE e^{\gl Y}$ as $n\to\infty$.
\end{proof}

\begin{corollary}
\label{C:simple}
For every $n\ge1$, we have
\begin{equation*}
\EE\,e^{\gl Y_n} \leq 
  \begin{cases} 
    e^{1.25\fudge^2 \gl^2},                    & \gl \leq 0, \\
    e^{0.5 \fudge^2\gl^2}, & -0.62\,n/(n+1) \leq \gl \leq 0,\\
 	e^{\fudge^2\gl^2}, & 0 \leq \gl \leq 0.42\,n/(n+1),\\
    e^{12\fudge^2 \gl^2},            & 0 \leq \gl \leq L_0\,n/(n+1),\\
    e^{2 e^{\fudge\gl}},               & L_0\,n/(n+1)\leq \gl.
   \end{cases}
\end{equation*}
In particular, 
$\EE\,e^{\gl Y_n} \leq 
\exp \left(\max \left( 12 \fudge^2\gl^2, 2 e^{\fudge\gl} \right)
\right)$ 
for all $\gl \in \RR$.
\end{corollary}
\begin{proof}
$\gl Y_n = \gl_n \Yh_n$ with $\gl_n:=\fudge\gl$.
\end{proof}

\begin{remark}\label{R:complicated}
The factors $\fudge$ in \refC{C:simple} are annoying but hardly
important in applications. With some effort, we have been able to
modify the proof in \cite{SJFill2} and obtain 
for $\gl\ge-0.58$ 
the same estimates
for $\EE e^{\gl Y_n}$ as obtained there for $\EE e^{\gl Y}$;
for $\gl<-0.58$ we only obtain a slightly weaker bound, which for large $n$
is inferior to the bound in \refC{C:simple}.
More precisely, we have shown
\begin{equation}
\EE\,e^{\gl Y_n} \leq 
  \begin{cases}
    e^{1.34 \gl^2},                 & \gl \leq -0.58, \\
    e^{0.5 \gl^2},                 & -0.58 \leq \gl \leq 0, \\
    e^{\gl^2},                     & 0 \leq \gl \leq 0.42, \\
    e^{12 \gl^2},                  & 0.42 \leq \gl \leq L_0, \\
    e^{2 e^{\gl}},                 & L_0 \leq \gl.
   \end{cases}
\end{equation}
In particular, $\EE\,e^{\gl Y_n} \leq \exp \left( \max \left( 12 \gl^2, 2
e^{\gl} \right) \right)$ for all $\gl \in \RR$.
In other words, we can
eliminate the factors $\fudge$ in \refC{C:simple} for $\gl\ge-0.58$
(and in particular for all positive $\gl$).
Since the proof is quite long and the result only
marginally improves \refC{C:simple}, we give the proof  not here
but rather in a separate appendix \cite{SJFill3App}.

It seems likely that with  further effort one could remove the factor
$\fudge$ for $\gl<-0.58$ too, so that all the bounds in \eqref{homer}
also would bound $e^{\gl Y_n}$. Moreover, it seems quite likely
that $\EE\,e^{\gl Y_n}\le \EE\,e^{\gl Y}$ holds for all $\gl$ and $n$,
and perhaps even that $\EE\,e^{\gl Y_n} \uparrow \EE\,e^{\gl Y}$,
as was proved for $\Yh_n$ in \refT{T:simple}.
\end{remark}

\refT{T:simple} enables us to get an explicit constant in
\Roesler's~\cite{Roesler} large deviation bound.

\begin{corollary}
\label{C:M1}
For
any $\eps > 0$ and $\gl > 0$,
$$
P(|X_n - \mu_n| \geq \eps \mu_n) \leq 2 \exp \left[ 3 \eps \gl + \max \left(
12 \gl^2, 2 e^{\gl} \right) \right] n^{- 2 \eps \gl}.
$$
\end{corollary}

\begin{proof}
By Markov's inequality,
\begin{align*}
P(|X_n - \mu_n| \geq \eps \mu_n)
  &  =  P(|\Yh_n| \geq \eps \mu_n / (n+1)) \\
  &\leq \exp(- \eps \gl \mu_n / (n + 1)) \EE\,e^{\gl |\Yh_n|} \\
  &\leq \exp(- \eps \gl \mu_n / (n + 1)) \left( \EE\,e^{\gl \Yh_n} + \EE\,e^{- \gl
           \Yh_n} \right).
\end{align*}
The result follows from \refT{T:simple}, since 
$\mu_n / (n+1) \geq 2 H_n - 4 \geq 2 \ln n - 3$
by~\eqref{castor} and~\eqref{jason}.  
\end{proof}

\begin{corollary}
For any fixed $\eps > 0$,
$$
P(|X_n - \mu_n| \geq \eps \mu_n) \leq n^{- 2 \eps \ln \ln n + O(1)}, \qquad n
\geq 2.
$$
\end{corollary}

\begin{proof}
Take (for $n\ge3$) $\gl = \ln \ln n$ in \refC{C:M1}.
\end{proof}

Finally we show that the rate of convergence of the moment generating
functions $\EE\,e^{\gl Y_n}$ to $\EE\,e^{\gl Y}$ also is $O(n^{- 1 / 2})$.
(The same holds for $\EE\,e^{\gl \Yh_n}$.)
\begin{theorem}
\label{T:mgfrate}
For any fixed complex~$\gl$,
$$
\EE\,e^{\gl Y_n} = \EE\,e^{\gl Y} + O(n^{- 1 / 2}).
$$
Explicitly, with $\gl_1 := \Re(\gl)$,
$$
\left| \EE\,e^{\gl Y_n} - \EE\,e^{\gl Y} \right| \leq 3 |\gl| \exp \left[ \max
\left( 24\fudge^2 \gl^2_1, e^{2 \fudge\gl_1} \right) \right] n^{- 1 / 2}.
$$
\end{theorem}

\begin{proof}
Consider a $d_2$-optimal coupling of~$Y_n$ and~$Y$.  Then, using the mean
value theorem, the \CS\ inequality, \refC{C:simple}, and~\eqref{homer},
\begin{eqnarray*}
\left| \EE\,e^{\gl Y_n} - \EE\,e^{\gl Y} \right|
  &\leq& \EE \left| e^{\gl Y_n} - e^{\gl Y} \right| \\
  &\leq& \EE \left( |\gl| |Y_n - Y|\,e^{\max(\gl_1 Y_n, \gl_1 Y)} \right) \\
  &\leq& |\gl| \left( \EE |Y_n - Y|^2  \right)^{1 / 2} \left( \EE\,e^{2
           \max(\gl_1 Y_n, \gl_1 Y)} \right)^{1 / 2} \\
  &\leq& |\gl| d_2(Y_n, Y) \left( \EE\,e^{2 \gl_1 Y_n} + \EE\,e^{2 \gl_1 Y}
           \right)^{1 / 2} \\
  &\leq& \sqrt{2}\,|\gl| \exp \left[ 
     \max \left( 24\fudge^2 \gl^2_1, e^{2\fudge \gl_1}
           \right) \right] d_2(Y_n, Y).
\end{eqnarray*}  
The result follws by \refT{T:d2}.
\end{proof}

\begin{remark}
By \refR{R:complicated}, the factors $\fudge$  can be eliminated in the
statement of \refT{T:mgfrate}. 
\end{remark}

{\bf Acknowledgment.\ }We thank Anhua Lin and Ludger Ruschendorf for helpful
discussions.

\end{document}